\documentclass{ACmod}

\DeclareFontFamily{U}{rsf}{}
\DeclareFontShape{U}{rsf}{m}{n}{
  <5> <6> rsfs5 <7> <8> <9> rsfs7 <10-> rsfs10}{}
\DeclareMathAlphabet{\mathscr}{U}{rsf}{m}{n}

\DeclareMathAlphabet{\mathgth}{U}{euf}{m}{n}

\DeclareFontFamily{U}{cyr}{}
\DeclareFontShape{U}{cyr}{m}{n}{
  <5> wncyr5 <6> wncyr6 <7> wncyr7 <8> wncyr8 <9> wncyr9 <10-> wncyr10}{}
\DeclareMathAlphabet{\mathcyr}{U}{cyr}{m}{n}

\input cyracc.def

\DeclareSymbolFont{bbold}{U}{bbold}{m}{n}
\DeclareSymbolFontAlphabet{\mathbbold}{bbold}

\makeatletter
\def\operator@font{\sf}
\makeatother

\newcommand{\cA}{{\mathscr A}}

\newcommand{\cF}{{\mathscr F}}

\newcommand{\cM}{{\mathscr M}}

\newcommand{\MF}{{\sf MF}}

\newcommand{\hres}{{\mathsf{hres}}}

\newcommand{\bone}{{\mathbf 1}}

\newcommand{\GM}{\mathsf{GGM}}

\DeclareMathOperator{\Tr}{Tr}

\DeclareMathOperator{\Spec}{Spec}

\DeclareMathOperator{\Der}{Der}

\DeclareMathOperator{\wt}{wt}

\DeclareMathOperator{\HC}{HC}

\DeclareMathOperator{\End}{End}

\DeclareMathOperator{\imag}{Im}

\DeclareMathOperator{\KS}{KS}

\newcommand{\ra}{\rightarrow}

\newcommand{\Q}{\mathbb{Q}}
\newcommand{\Z}{\mathbb{Z}}

\renewcommand{\phi}{\varphi}

\usepackage{tikz}
\usepackage{amscd}
\usepackage[tocflat]{tocstyle}
\usetikzlibrary{positioning}
\numberwithin{equation}{section}

\begin{document}

\title{Categorical primitive forms and Gromov-Witten invariants of $A_n$ singularities}

\author[C\u ald\u araru, Li and Tu]{%
Andrei C\u ald\u araru\footnote{Partially supported by the National
Science Foundation through grant numbers DMS-1200721 and DMS-1811925.}\;\; Si Li \footnote{Partially supported by Grant 20151080445 of Independent Research
Program at Tsinghua University}\;\; Junwu Tu\footnote{Partially supported by the National Science Foundation through grant number DMS-1801806.}}

\address{Andrei C\u ald\u araru, Mathematics Department,
University of Wisconsin--Madison, 480 Lincoln Drive, Madison, WI
53706--1388, USA.}
\address{Junwu Tu, Institute of Mathematical Science, ShanghaiTech University, Shanghai 201210, China.} 
\address{Si Li, Mathematics Department and Yau Mathematical Sciences Center, Tsinghua University, Beijing 100084, China.}

\begin{abstract} 
  {\sc Abstract:} We introduce a categorical analogue of Saito's notion of primitive forms. For the category $\MF(\frac{1}{n+1}x^{n+1})$ of matrix factorizations of $\frac{1}{n+1}x^{n+1}$, we prove that there exists a unique, up to non-zero constant, categorical primitive form. The corresponding genus zero categorical Gromov-Witten invariants of $\MF(\frac{1}{n+1}x^{n+1})$ are shown to match with the invariants defined through unfolding of singularities of $\frac{1}{n+1}x^{n+1}$.
\end{abstract}

\maketitle

\setcounter{tocdepth}{2}

\section{Introduction}

\paragraph{{\bf Classical versus homological mirror symmetry.}}
In the seminal 1991 paper~\cite{COGP} the physicists Candelas-de la
Ossa-Green-Parkes proposed a mysterious procedure for calculating the
$g=0$ Gromov-Witten invariants of the quintic Calabi-Yau threefold.
Their procedure matched these invariants with the coefficients
of a solution of a certain differential equation governing the
variation of Hodge structures of a mirror family of Calabi-Yau
threefolds. There was no mathematical explanation for the validity of
this correspondence, as the motivation for it came from string theory.

Three years later, in his ICM talk~\cite{Kon}, Kontsevich proposed a
conjectural explanation (the {\em homological mirror symmetry
  conjecture}) for the physicists' computation. He envisioned the
existence of an equivalence between certain triangulated categories
defined on the two sides of the mirror. He further suggested that the
equality of invariants predicted by physics should follow from the
existence of certain Gromov-Witten-like invariants associated to these
categories.

\paragraph{{\bf The VSHS versus the TCFT approaches to categorical invariants.}}
Two approaches for the construction of such invariants have emerged
since Kontsevich's original proposal. The first of these, which is
only useful for recovering genus zero invariants, involves the study
of a certain Hodge structure first introduced by
Saito~\cite{Sai2}~\cite{Sai1} in his study of the theory of primitive
forms.  Later Barannikov~\cite{Bar},~\cite{Bar2} and
Barannikov-Kontsevich~\cite{BarKon}) introduced the important notion
of (polarized) Variation of Semi-infinite Hodge Structures (VSHS),
generalizing Saito's framework to many other geometric contexts.  These
works set up the stage for the study of the analogous
non-commutative/categorical Hodge structures that will be studied in
this paper, see also~\cite{GPS} and~\cite{She}.

\[\begin{tikzpicture}[box/.style={draw,rounded corners,text width=2.5cm,align=center}]
\node[box] at (-5,0) (a) {Compact, smooth CY $A_\infty$ categories};
\node[box] at (0,0) (b) {VSHS's};
\node[box] at (5,0) (d) {Genus zero categorical Gromov-Witten theory};
\draw[thick,->] (a) -- (b) node[midway,above]{A};
\draw[thick,->] (b) -- (d) node[midway,above]{B};
\end{tikzpicture}\] Step A in this diagram is well understood due to
the works of Katzarkov-Kontsevich-Pantev~\cite{KKP},
Shklyarov~\cite{Shk}~\cite{Shk2}, Sheridan~\cite{She}.
Ganatra-Perutz-Sheridan~\cite{GPS} obtained a partial resolution of
Step B for a particular type of VSHS's called {{\sl Hodge-Tate VSHS's
    over a one dimensional deformation space}}. Remarkably, this
already suffices for Ganatra-Perutz-Sheridan to give a conceptual
explanation of the computation by Candelas-de la Ossa-Green-Parkes
using homological mirror symmetry for quintics (another important
result of Sheridan~\cite{She2}).

A second approach, which works
for invariants of arbitrary genus, is due to
Kontsevich-Soibelman~\cite{KonSoi} and
Costello~\cite{Cos1}~\cite{Cos}.  It proceeds by constructing a
2-dimensional TCFT from a given Calabi-Yau $A_\infty$ category, and
extracting Gromov-Witten invariants from this TCFT.
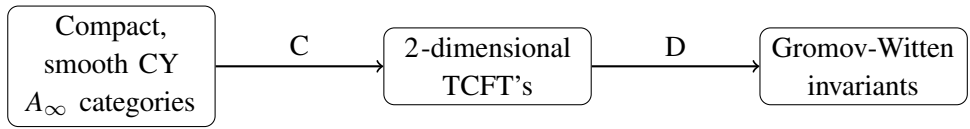
\[\begin{tikzpicture}[box/.style={draw,rounded corners,text width=2.5cm,align=center}]
\node[box] at (-5,0) (a) {Compact, smooth CY $A_\infty$ categories};
\node[box] at (0,0) (c) {$2$-dimensional TCFT's};
\node[box] at (5,0) (d) {Gromov-Witten invariants};
\draw[thick, ->] (a) -- (c) node[midway,above]{C};
\draw[thick,->] (c) -- (d) node[midway,above]{D};
\end{tikzpicture}\] Unfortunately, aside from the foundational works
mentioned above, little progress has been made in understanding these
invariants. One reason is the fact that explicit computation in this
approach is extremely difficult. At the time of the writing of this
paper the only computed case is the $g=1$, $n=1$ categorical invariant
of the derived category $D^b_{\sf{coh}} (E_\tau)$ of an elliptic curve
$E_\tau$, obtained by the first and third authors~\cite{CalTu}.

\paragraph{{\bf The present paper.}} 
In this paper we initiate a first case study of categorical
Gromov-Witten invariants (in both the VSHS and TCFT approaches) for
categories of matrix factorizations. We shall focus on the category
$\MF(\frac{1}{n+1} x^{n+1})$ of matrix factorizations of
$\frac{1}{n+1} x^{n+1}$. (The coefficient $\frac{1}{n+1}$ is a
normalization to make later computations easier.) 

\paragraph{{\bf Primitive forms.}} 
We begin with the VSHS approach. In the case of categories of matrix
factorizations the method of~\cite{GPS} does not apply because the
deformation space we shall work with is not one dimensional, and we
are not in the Hodge-Tate case. We may obtain useful insights from the
original works of K.\ Saito~\cite{Sai2},~\cite{Sai1}. He gave a
construction of genus zero Gromov-Witten type invariants associated
with a holomorphic function $f(x_1,\cdots,x_n)$ with an isolated
singularity at $x_1=\cdots=x_n=0$. A key ingredient of Saito's theory
are the so-called primitive forms. The definition of primitive forms
is quite involved: in particular, their existence is a highly
non-trivial question. To prove the existence it was observed by K.\
Saito~\cite{Sai1} that there exists a one-to-one correspondence
between the set of primitive forms (up to a constant rescaling) and
the set of {{\sl good sections}}, or equivalently {{\sl good
    splittings}} of the Hodge filtration on the vanishing cohomology
of the function $f(x_1,\ldots,x_n)$. This allowed him to construct
primitive forms for weighted homogeneous polynomials. The existence of
primitive forms for arbitrary isolated singularity was proved by M.\
Saito~\cite{SaiM} using mixed Hodge theory.

Explicit expressions for primitive forms are rarely known, except in a handful of examples: ADE and simple elliptic singularities~\cite{Sai1}.  There is also a formal version of primitive forms developed by Li-Li-Saito~\cite{LLS}.  A recursive algorithm ~\cite{LLS}~\cite{LLSS} was obtained to explicitly compute primitive forms perturbatively, culminating with a proof of Landau-Ginzburg mirror symmetry for exceptional unimodular singularities~\cite{LLSS} and for general invertible polynomials~\cite{HSLW}.

\paragraph{{\bf Categorical primitive forms.}} Motivated by the previous discussion, for a Calabi-Yau $A_\infty$ category we introduce 
\begin{itemize}
\item[(1)] a categorical analogue of Saito's notion of primitive forms (Definition~\ref{def:primitive});
\item[(2)] the notion of a {{\sl good splitting}} of the Hodge filtration (Definition~\ref{def:goodsplitting}) on the periodic cyclic homology.
\end{itemize} 
An important ingredient in these definitions is a formula (Equation~\eqref{eq:u-connection}) for the canonical $u$-direction connection on the periodic cyclic homology of a $\Z/2\Z$-graded $A_\infty$ algebra. In the differential graded case, this $u$-connection was conceptually described by Katzarkov-Kontsevich-Pantev~\cite{KKP}, and an explicit formula was written down later by Shklyarov~\cite{Shk}. 

We remark that the notion of good splitting in $(2)$ is different from Costello's original definition~\cite[Definition 11.0.11]{Cos} of a splitting of the Hodge filtration. In the $\Z$-graded (or even $\Q$-graded) case, our notion is equivalent to Costello's definition with an additional requirement of homogeneity. In the $\Z/2\Z$-graded case, the condition for a splitting to be good is more mysterious: it sort of requires maximal homogeneity even if we do not have a $\Z$-grading (or $\Q$-grading) to begin with. Furthermore, a splitting must also be compatible with the Calabi-Yau structure. A precise compatibility condition can be found in Definition~\ref{def:goodsplitting}.
 
As a first example we study categorical primitive forms for the category $\MF(\frac{1}{n+1} x^{n+1})$. We prove that there exists a unique categorical primitive form up to non-zero scalars. We then use the aforementioned Li-Li-Saito algorithm to compute the corresponding genus zero categorical invariants. As one may expect, the categorical computation matches with that of Saito's theory of unfolding of singularities of $\frac{1}{n+1} x^{n+1}$. See Theorem~\ref{thm:comparison} for details.

\paragraph{{\bf Towards $B$-model Gromov-Witten invariants of Landau-Ginzburg orbifolds.}}
Let $(A,\omega)$ be a smooth and proper Calabi-Yau $A_\infty$ algebra. In a forthcoming work~\cite{CalTu2} the first and the third named authors obtain a bijection 
\[
\begin{tikzpicture}[box/.style={draw,rounded corners,text width=4cm,align=center}]
\node[box] at (-3,0) (a) {Categorical primitive forms of $(A,\omega)$};
\node[box] at (3,0) (d) {$\omega$-compatible good splittings of $HP_*(A)$ };
\draw [thick,->] (a) -- (d) node[midway,above]{$\cong$} ;
\end{tikzpicture}
\]
This result is the categorical analogue of the aforementioned bijections of K.\ Saito, M.\ Saito and Li-Li-Saito. It can be used to obtain a construction of the genus zero $B$-model theory for the category $\MF(W,G)$ of $G$-equivariant matrix factorizations. Conjecturally, this theory should be mirror dual to the genus zero part of the Fan-Jarvis-Ruan-Witten (FJRW) theory~\cite{FJRW}.

\paragraph{{\bf TCFT invariants.}}
We also investigate in
Propositions~\ref{prop:vanishing},~\ref{prop:11}, and~\ref{prop:03}
Costello's $(0,3)$ and $(1,1)$ invariants for the category
$\MF(\frac{1}{n+1} x^{n+1})$. The results match the expected ones from
FJRW theory. For higher genus invariants it remains a difficult task
to perform explicit computations in Costello's categorical approach.

An important ingredient in Costello's definition of categorical Gromov-Witten invariants is a choice of a splitting of the Hodge filtration on periodic cyclic homology.  Such a choice is required to satisfy certain properties, see Costello~\cite{Cos}.  We propose the following conjecture, which we have verified for the category $\MF(\frac{1}{n+1} x^{n+1})$.
\vspace{-0.3em}

\begin{Conjecture} Let $(A, \omega)$ be a smooth and proper Calabi-Yau
  category. Under the bijection above Costello's genus zero categorical invariants of $(A,\omega)$ match with the invariants defined using the corresponding primitive form.  
\end{Conjecture}
\medskip

\paragraph{{\bf Contents of the paper.}}
The paper is organized as follows. In Section~\ref{sec:invariants} we identify an explicit splitting of the Hodge filtration of the category $\MF(\frac{1}{n+1} x^{n+1})$. This allows us to compute Costello's $(0,3)$, and $(1,1)$ invariants of this category. In Section~\ref{sec:primitive} we write down a formula for the canonical $u$-direction connection on the periodic cyclic homology of a $\Z/2\Z$-graded $A_\infty$ algebra. Then we introduce a notion of categorical primitive forms, and prove their existence and uniqueness for the category $\MF(\frac{1}{n+1} x^{n+1})$. In Section~\ref{sec:computation} we match the genus zero categorical Gromov-Witten invariants associated to the canonical primitive form of $\MF(\frac{1}{n+1} x^{n+1})$ defined in the previous section, with Saito's invariants obtained through unfolding of singularities.

\paragraph{{\bf Notations and Conventions.}} We work with the shifted sign conventions for $A_\infty$ algebras: the structure maps $\mu_k$'s are of degree $1$ for all $k\geq 0$. On the other hand, for the Hochschild chain complex we always use homological grading, i.e. $b$ has degree $-1$.

\paragraph {\bf Acknowledgments.} We thank Nick Addington, Kevin Costello and Sasha Polishchuk for useful discussions. 
\section{Invariants of $(A_n,\omega)$}
\label{sec:invariants}

In this section, we explicitly construct a splitting of the Hodge
filtration of the periodic cyclic homology of
$\MF(\frac{1}{n+1}x^{n+1})$. This allows to compute the $(0,3), (1,1)$
invariants of this category invariants in Costello's framework.

\paragraph{{\bf The minimal model of $\MF(\frac{1}{n+1}x^{n+1})$.}}
Consider the potential function
$W: \mathbb{A}^1 \rightarrow \mathbb{A}^1$ given by
$W(x)=\frac{1}{n+1}x^{n+1}$. A result of Dyckerhoff~\cite{Dyc} shows
that its category of matrix factorizations is compactly generated by
the stabilization of the residue field
\[ \mathbb{K}^{{\sf stab}}:=\mathbb{K}[x]
  \stackrel{x}{\longrightarrow} \mathbb{K}[x] \stackrel{\frac{1}{n+1}
    x^n}{\longrightarrow} \mathbb{K}[x].\] 
The differential graded algebra $\End(\mathbb{K}^{{\sf stab}})$ can be
described explicitly as follows:
\[ \End(\mathbb{K}^{{\sf stab}})=\mathbb{K}[x]\otimes_{\mathbb{K}}
  \mathbb{K}\langle \hat{x}, \hat{y}\rangle, \] where
$\hat{x}, \hat{y}$ are odd variables satisfying the commutator
relation $\hat{x}\hat{y}+\hat{y}\hat{x}=1$. The differential in this
dg algebra is $\mathbb{K}[x]$-linear and acts on generators by $d(\hat{x})=x$
and $d(\hat{y})=\frac{1}{n+1} x^n$.

As described in~\cite{Dyc} one can apply the homological perturbation
lemma to the above dg algebra in order to obtain its minimal model. This gives
an $A_\infty$ algebra which we shall denote by $A_n$.  Explicitly, the
algebra $A_n={{\sf span}}( \bone, \epsilon)$ is a $2$-dimensional
$\Z/2\Z$-graded vector space generated by $\bone$ even, and $\epsilon$
odd. The $A_\infty$ multiplications are determined by the following
properties:
\begin{itemize}
\item[(a)] $\bone$ is a strict unit,
\item[(b)] apart from (a), the only non-zero multiplication is given by
\[ \mu_{n+1}(\epsilon, \ldots, \epsilon)=\frac{1}{n+1}\bone.\]
\end{itemize}
We also put an odd inner product on $A_n$ defined by $\langle \bone, \epsilon\rangle=1$. One verifies that this pairing makes the functional
\[ a_1\otimes \cdots a_{k+1} \mapsto \langle \mu_k(a_1,\cdots,a_k), a_{k+1}\rangle\]
invariant under cyclic permutations.  Thus we obtain a cyclic $A_\infty$ algebra structure. This cyclic structure is natural in the sense that any odd inner product is a constant multiple of $\langle-,-\rangle$.

\paragraph{{\bf Hochschild invariants.}}
The projection $\pi: A_n \ra \mathbb{K}$ defined by
$\pi(a\bone+b\epsilon)=a$ is an augmentation of $A_n$. Using this
augmentation we form the reduced Hochschild chain complex:
\[ \cdots \longrightarrow \bone|\epsilon|\cdots|\epsilon \oplus \epsilon|\epsilon|\cdots|\epsilon \longrightarrow \cdots \bone|\epsilon|\epsilon\oplus \epsilon|\epsilon|\epsilon \longrightarrow \bone|\epsilon\oplus \epsilon|\epsilon\longrightarrow \bone \oplus \epsilon.\]
For convenience, let us denote by $\bone|\epsilon^k$ the element 
$$\bone|\underbrace{\epsilon|\cdots|\epsilon}_{\text{$k$-copies}}.$$
The parity of $\bone|\epsilon^k$ in the Hochschild complex is even for
all $k\geq 0$. Similarly define $\epsilon|\epsilon^k$; it has odd
parity.

The Hochschild differential is given by
\begin{align*}
 b & \big( \bone|\epsilon^k\big) = 0\\
 b & \big( \epsilon|\epsilon^k\big)= \begin{cases}
 0, \mbox{ \; if \;} k< n\\
 \bone|\epsilon^{k-n}, \mbox{\; if \; } k\geq n\\
 \end{cases}
 \end{align*}
These formulas immediately imply the following proposition.

\begin{Proposition}~\label{prop:homology}
The Hochschild homology of $A_n$ is an $n$-dimensional vector space over $\mathbb{K}$, concentrated in odd degree. Furthermore, we have the following set of explicit generators
\[ HH_*(A_n)= {{\sf span}}\big( [\epsilon], [\epsilon|\epsilon], \cdots, [\epsilon|\epsilon^{n-1}]\big).\]
\end{Proposition}

\paragraph{{\bf The Frobenius algebra structure.}}
Kontsevich-Soibelman~\cite{KonSoi} and Cos\-te\-llo~\cite{Cos1} exhibited
an important structure on the shifted Hochschild complex
$CC_{*+d}(A)$: it carries an action of the PROP of (a shifted version
of) chains of moduli spaces of Riemann surfaces with parametrized
boundaries. Namely, we have a chain map
\[ \rho : C^{{\sf comb}}_{*+d(2-2g-2n)}(\cM_{g,n,m})\otimes
  CC_{*+d}(A)^{\otimes n} \rightarrow CC_{*+d}(A)^{\otimes m}.\] Here,
the chain complex $C^{{\sf comb}}_*(\cM_{g,n,m})$ with $n>0$ denotes a
certain combinatorial model of $\cM_{g,n,m}$, the moduli space of
curves of genus $g$ with $n$ parametrized input boundaries and $m$
parametrized output boundaries.  This combinatorial model uses
decorated ribbon graphs, and calculates the rational homology of the
moduli space $\cM_{g,n,m}$. 

In particular we have a pairing of degree $2d$ defined by
\[\rho\big( \begin{tikzpicture}[baseline={(current bounding box.center)},scale=0.2]
\draw [thick] (0,2) circle [radius=2];
\draw [thick] (-2,2) to (-0.6,2);
\draw [thick] (-0.8,2.2) to (-0.4,1.8);
\draw [thick] (-0.8, 1.8) to (-0.4, 2.2);
\draw [thick] (2,2) to (3.4,2);
\draw [thick] (3.2,2.2) to (3.6,1.8);
\draw [thick] (3.2, 1.8) to (3.6, 2.2);
\end{tikzpicture}  \big) : CC_{*+d}(A)^{\otimes 2} \rightarrow \mathbb{K}.\]
This is known as the chain level Mukai pairing. We shall denote this
map by $\langle-,-\rangle_{{\sf Muk}}$. One can show this pairing is
always symmetric. The fact that this agrees with the categorical
definition of Shklyarov~\cite{Shk} is a result of
Sheridan~\cite{She}.

\begin{Proposition}~\label{prop:pairing}
Let $CC_{*+1}(A_n)$ be the shifted Hochschild chain complex of the
$A_\infty$ algebra $A_n$. Then the chain level  Mukai pairing is given by
\[ \langle \epsilon|\epsilon^i, \epsilon|\epsilon^j\rangle_{{\sf Muk}}=\begin{cases}
1, \mbox{\; if \;}  i+j=n-1\\
0, \mbox{\; otherwise.}
\end{cases}\]
In particular, it induces a symmetric non-degenerate pairing on the purely even vector space $HH_{*+1}(A_n)= {{\sf span}}\big( [\epsilon], [\epsilon|\epsilon], \cdots, [\epsilon|\epsilon^{n-1}]\big)$.
\end{Proposition}

\begin{proof}
The insertions $\epsilon^i$ and $\epsilon^j$ must all go to the same
vertex in order to produce a nonzero evaluation. Since only $m_{n+1}$ is not zero, we need to have $i+j=n-1$. There are $i+1$ and $j+1$ ways to arrange the insertions to each of the two vertices, and each produces an evaluation of $\frac{1}{n+1}$. Thus, if $i+j=n-1$,  we have
\[  \langle \epsilon[\epsilon^i], \epsilon[\epsilon^j]\rangle_{{\sf Muk}}= \frac{1}{n+1}\big( i+1+j+1\big)=1.\]
\end{proof}

One can also compute the coproduct of $CC_{*+1}(A_n)$ (dual to the cup product of $CC^*(A_n)$) using the graph
\[\rho\big(\begin{tikzpicture}[baseline={([yshift=-2ex]current bounding box.center)},scale=0.2]
\draw [thick] (0,0) to (0,2);
\draw [thick] (-0.2, 1.8) to (0.2, 2.2);
\draw [thick] (0.2, 1.8) to (-0.2, 2.2);
\draw [thick] (0,0) to (-2,0);
\draw [thick] (0,0) to (2,0);
\draw [thick] (-2.2,0) circle [radius=0.2];
\draw [thick] (2.2,0) circle [radius=0.2];
\end{tikzpicture}\big): CC_{*+1}(A_n)\ra CC_{*+1}(A_n)^{\otimes 2}\]

\begin{Proposition}~\label{prop:coproduct}
The coproduct of $H$ is given by formula
\[ \rho\big(\begin{tikzpicture}[baseline={([yshift=-2ex]current bounding box.center)},scale=0.2]
\draw [thick] (0,0) to (0,2);
\draw [thick] (-0.2, 1.8) to (0.2, 2.2);
\draw [thick] (0.2, 1.8) to (-0.2, 2.2);
\draw [thick] (0,0) to (-2,0);
\draw [thick] (0,0) to (2,0);
\draw [thick] (-2.2,0) circle [radius=0.2];
\draw [thick] (2.2,0) circle [radius=0.2];
\end{tikzpicture}\big) \big( [\epsilon|\epsilon^k]\big) = \sum_{i+j=k} [\epsilon|\epsilon^i]\otimes [\epsilon|\epsilon^j ] .\]
\end{Proposition}

\begin{proof}
The sum comes from inserting $k$ copies of $\epsilon$'s into the two
outgoing white vertices.
\end{proof}

\paragraph{{\bf Splittings of the Hodge filtration.}} 
The chain level Mukai-pairing induces the so-called higher residue
pairing on the periodic cyclic chain complex:
 \begin{align*}
 \langle -,-\rangle_\hres & :~  CC_{*+d}(A)\laurent{u}^{\otimes 2} \rightarrow \mathbb{K}\laurent{u},\\
\langle \alpha u^i, \beta u^j \rangle_\hres &:= (-1)^i \langle
  \alpha, \beta \rangle_{{\sf Muk}} \cdot u^{i+j}, \;\; \forall \alpha, \beta\in CC_{*+d}(A), \;\; i, j\in \Z.
\end{align*}
The higher residue pairing is a chain map with respect to the cyclic
differential $b+uB$, which makes it descend to a pairing on periodic
cyclic homology. The space $CC_{*+d}(A)\laurent{u}$ is filtered by powers of $u$, and
the pairing respects this filtration.  In particular, for
$x, y$ in the negative cyclic homology $\HC_{*+d}^-(A)$ of $A$ we
have $\langle x, y \rangle_\hres \in \mathbb{K}\series{u}$.

The following definition was introduced in~\cite{CalTu}.

\begin{Definition}~\label{defi-splitting}
Denote by $ HH_{*+d}(A):= H_*\big(CC_{*+d}(A)\big)$ the shifted Hochschild homology of $A$. A splitting of the Hodge filtration on $HP_*(A)$ is a linear map
\[ s: HH_{*+d}(A) \rightarrow HC^-_{*+d}(A),\]
such that 
\begin{itemize}
\item[(1)] it splits the canonical (constant term) map $HC_{*+d}^-(A) \rightarrow HH_{*+d}(A)$,
\item[(2)] for any two classes $\alpha,\beta\in HH_{*+d}(A)$, we have $\langle \alpha,\beta \rangle_{{\sf Muk}} = \langle s(\alpha), s(\beta)\rangle_\hres$.
\end{itemize}
\end{Definition}
\medskip
 
The second condition is equivalent to requiring that the pairing
$\langle s(\alpha), s(\beta)\rangle_\hres$, which {{\sl a priori}}
takes values in $\mathbb{K}\series{u}$, be a constant in $\mathbb{K}$.

\paragraph{{\bf Explicit splitting of $HP_*(A_n)$.}}
To compute the categorical Gromov-Witten invariants of $A_n$ we need
to split the Hodge filtration. Since the Hochschild homology
$HH_{*+1}(A_n)$ is concentrated in the even degree, the
Hodge-de Rham spectral sequence always degenerates. However, the
choice of a splitting is not unique. To fix this ambiguity we put
rational weights on $\bone$ and $\epsilon$. The existence of this
additional structure is due to the homogeneity of the potential
function $W=\frac{1}{n+1}x^{n+1}$.

More precisely, we put the following fractional weight on $A_n$:
\[ \wt(\epsilon)=\frac{n-1}{n+1}, \;\;\wt(\bone)=0. \]
With respect to these weights the multiplication $m_k$ on $A_n$ is of weight $2-k$.  This induces the corresponding weight on the Hochschild complex:
\[ \wt( \bone|\epsilon^k) = \frac{2k}{n+1},\;\; \wt(\epsilon| \epsilon^k)= \frac{2k+1-n}{n+1}.\]
The Hochschild differential $b$ is of weight $-1$, and the cyclic differential $B$ is of weight $1$.
We also set $\wt(u)=-2$ for the circle parameter. This makes the
cyclic differential $b+uB$ of total weight $-1$.

\begin{Proposition}~\label{prop:splitting} For the $A_\infty$ algebra
  $A_n$ there exists a unique weight preserving splitting of the Hodge
  filtration. Explicitly, it is given by the formula
\[ s([\epsilon|\epsilon^k])= \sum_{l\geq 0}^\infty (-1)^l c_{k,l} \;\epsilon|\epsilon^{k+nl+l} u^l,\]
where for each $0\leq k\leq n-1$ we set $c_{k,0}=1$ and $c_{k,l}:=\prod_{j=0}^{l-1} \big((k+1+j(n+1)\big)$ for any $l\geq 1$.
\end{Proposition}

\begin{proof}
One can directly check that the above splitting is
weight-preserving. Furthermore, at each stage of lifting from working
modulo $u^l$
to modulo $u^{l+1}$ the ambiguity is an element of $HH_*(A_n)= {{\sf
    span}}\big( [\epsilon], [\epsilon|\epsilon], \cdots,
[\epsilon|\epsilon^{n-1}]\big)$. But the weights are all inside the
interval $[-\frac{n-1}{n+1},\frac{n-1}{n+1}]$. 
Thus requiring the splitting to be weight preserving kills the
ambiguity by degree reason.
\end{proof}

\paragraph{{\bf Some computations.}}
We use the above splitting to compute Costello's categorical
Gromov-Witten invariants of type $(0,3)$ and $(1,1)$. We start with
the $(1,1)$ invariants.  

\begin{Proposition}~\label{prop:vanishing}
For any $0\leq k\leq n-1$, we have the type $(1,1)$ categorical Gromov-Witten invariant $\langle[\epsilon|\epsilon^k]\rangle_{1,1}^{A_n}=0$.
\end{Proposition}
\begin{proof}
  (For general details on Costello's formalism see~\cite{Cos}
  and~\cite{CalTu}.) The computation takes place in the deformed Fock
  module. Denote by
\[ \widetilde{\epsilon|\epsilon^k}:= \sum_{l\geq 0}^\infty (-1)^l c_{k,l} \;\epsilon|\epsilon^{k+nl+l} u^l\]
the weight preserving lift, considered as an element in the periodic cyclic chain complex of $A_n$. The deformed Fock module relation gives
\[ \widetilde{\epsilon|\epsilon^k}  =  \left ( T_1(\epsilon|\epsilon^k) +
    T_2\big( -(k+1)\epsilon|\epsilon^{k+n+1}\big) + T_3(\epsilon|\epsilon^k) \right )\lambda +
  O(\lambda^2), \]
where the three terms are given by
\begin{eqnarray*}
 T_1(\epsilon|\epsilon^k) =&
\dfrac{1}{2}  \ \ \rho\big(
              \begin{tikzpicture}[baseline={([yshift=-2ex]current
                  bounding box.center)},scale=0.3] 
\draw [thick] (0,0) to (0,2);
\draw [thick] (-0.2, 1.8) to (0.2, 2.2);
\draw [thick] (0.2, 1.8) to (-0.2, 2.2);
\draw [thick] (0,0) to (-2,0);
\draw [thick] (0,0) to (2,0);
\draw [thick] (-2.2,0) circle [radius=0.2];
\draw [thick] (2.2,0) circle [radius=0.2];
\end{tikzpicture}\big) (\epsilon|\epsilon^k) \\\ \\
 T_2\big( -(k+1)\epsilon|\epsilon^{k+n+1}\big) =&
\dfrac{1}{24} \ \rho \big(\begin{tikzpicture}[baseline={(current
                  bounding box.center)},scale=0.3] 
\draw [thick] (0,2) circle [radius=2];
\draw [thick] (-2,2) to (-0.6,2);
\draw [thick] (-0.8,2.2) to (-0.4,1.8);
\draw [thick] (-0.8, 1.8) to (-0.4, 2.2);
\draw [thick] (-1.4142, 0.5858) to [out=40, in=140] (1, 0.5);
\draw [thick] (1.25, 0.3) to [out=-45, in=225] (2, 0.2);
\draw [thick] (2,0.2) to [out=45, in=-50] (1.732, 1);
\end{tikzpicture}\big) ( -(k+1)\epsilon|\epsilon^{k+n+1})\\\ \\
T_3(\epsilon|\epsilon^k)  =& 
- \dfrac{1}{4}\ \rho\big( 
\begin{tikzpicture}[baseline={(current bounding box.center)},scale=0.3]
\draw [thick] (0,2) circle [radius=2];
\draw [thick] (0,0) to (0,1.4);
\draw [thick] (-0.2, 1.2) to (0.2, 1.6);
\draw [thick] (-0.2, 1.6) to (0.2, 1.2);
\draw [thick] (0,0) to [out=80, in=180] (0.5, 1);
\draw [thick] (0.5,1) to [out=0, in=100] (0.9, 0.4);
\draw [thick] (0,0) to [out=-80, in=180] (0.5, -1);
\draw [thick] (0.5, -1) to [out=0, in=-100] (0.9, 0);
\end{tikzpicture}\big) (\epsilon|\epsilon^k).
\end{eqnarray*}
One verifies that the second and the third terms both give zero. The
first term is computed in Proposition~\ref{prop:coproduct}. Putting
everything together we get 

\begin{equation}\label{eq:1-insertion}
\begin{split}
\widetilde{\epsilon|\epsilon^k} &= \Big(\frac{1}{2} \sum_{i+j=k} \epsilon|\epsilon^i u^{-1} \cdot \epsilon|\epsilon^j u^{-1}\Big) \lambda + O(\lambda^2)\\
&= \Big(\frac{1}{2} \sum_{i+j=k} (u^{-1}\widetilde{\epsilon|\epsilon^i} + (i+1)\epsilon|\epsilon^{i+n+1}+\cdots) \cdot u^{-1}\widetilde{\epsilon|\epsilon^j}\Big) \lambda + O(\lambda^2)\\
&= \frac{1}{2} \sum_{i+j=k} (u^{-1}\widetilde{\epsilon|\epsilon^i})(u^{-1}\widetilde{\epsilon|\epsilon^j})\lambda+ \sum_{i+j=k} \frac{1}{2}\langle (i+1)\epsilon|\epsilon^{i+n+1}, \epsilon|\epsilon^j\rangle_{{\sf Muk}}\lambda+ O(\lambda^2)\\
&= \frac{1}{2} \sum_{i+j=k} (u^{-1}\widetilde{\epsilon|\epsilon^i})(u^{-1}\widetilde{\epsilon|\epsilon^j})\lambda +O(\lambda^2).
\end{split}
\end{equation}
The last equality follows from Proposition~\ref{prop:pairing}. The
invariant $\langle [\epsilon|\epsilon^k]\rangle_{1,1}$ is given by the
constant term coefficient of $\lambda$ of
$\widetilde{\epsilon|\epsilon^k}$ in the deformed Fock module. The
above computation shows that there is no constant term in the
expansion in powers of $\lambda$.
\end{proof}

\begin{Proposition}~\label{prop:11}
The type $(1,1)$ categorical Gromov-Witten invariant of $A_n$ with
gravitational descendant $u^l$ is given by
\[ \langle [\epsilon|\epsilon^k] u^l \rangle_{1,1}^{A_n}=\begin{cases}
\frac{n}{24}, \mbox{\; if \;} k=n-1, l=1;\\
0, \mbox{\;\; otherwise.}
\end{cases}\]
\end{Proposition}

\begin{proof}
The result comes from the evaluation of 
\[T_2\big( \epsilon|\epsilon^{k}\big) =
\dfrac{1}{24} \ \rho \big(\begin{tikzpicture}[baseline={(current
                  bounding box.center)},scale=0.3] 
\draw [thick] (0,2) circle [radius=2];
\draw [thick] (-2,2) to (-0.6,2);
\draw [thick] (-0.8,2.2) to (-0.4,1.8);
\draw [thick] (-0.8, 1.8) to (-0.4, 2.2);
\draw [thick] (-1.4142, 0.5858) to [out=40, in=140] (1, 0.5);
\draw [thick] (1.25, 0.3) to [out=-45, in=225] (2, 0.2);
\draw [thick] (2,0.2) to [out=45, in=-50] (1.732, 1);
\end{tikzpicture}\big) ( \epsilon|\epsilon^{k})\]
But this action can be decomposed into 
\[ \rho\big( \begin{tikzpicture}[baseline={(current
                  bounding box.center)},scale=0.3] 
\draw [thick] (0,2) circle [radius=2];
\draw [thick] (-2,2) to (-0.6,2);
\draw [thick] (-0.8,2.2) to (-0.4,1.8);
\draw [thick] (-0.8, 1.8) to (-0.4, 2.2);
\draw [thick] (-1.4142, 0.5858) to [out=40, in=140] (1, 0.5);
\draw [thick] (1.25, 0.3) to [out=-45, in=225] (2, 0.2);
\draw [thick] (2,0.2) to [out=45, in=-50] (1.732, 1);
\end{tikzpicture}\big) = \rho\big(\begin{tikzpicture}[baseline={(current bounding box.center)},scale=0.3]
\draw [thick] (0,2) circle [radius=2];
\draw [thick] (-2,2) to (-0.6,2);
\draw [thick] (-0.8,2.2) to (-0.4,1.8);
\draw [thick] (-0.8, 1.8) to (-0.4, 2.2);
\draw [thick] (2,2) to (3.4,2);
\draw [thick] (3.2,2.2) to (3.6,1.8);
\draw [thick] (3.2, 1.8) to (3.6, 2.2);
\end{tikzpicture}  \big) \rho \big( \begin{tikzpicture}[baseline={([yshift=-2ex]current bounding box.center)},scale=0.3]
\draw [thick] (0,0) to (0,2);
\draw [thick] (-0.2, 1.8) to (0.2, 2.2);
\draw [thick] (0.2, 1.8) to (-0.2, 2.2);
\draw [thick] (0,0) to (-2,0);
\draw [thick] (0,0) to (2,0);
\draw [thick] (-2.2,0) circle [radius=0.2];
\draw [thick] (2.2,0) circle [radius=0.2];
\end{tikzpicture}\big)\]
The computation now follows from Propositions~\ref{prop:pairing} and~\ref{prop:coproduct}.
\end{proof}

\begin{Proposition}~\label{prop:03}
The type $(0,3)$ categorical Gromov-Witten invariant of $A_n$ is given by formula
\[ \langle [\epsilon|\epsilon^i], [\epsilon|\epsilon^j], [\epsilon|\epsilon^k]\rangle^{A_n}_{0,3}=\begin{cases}
1, \mbox{\;\; if\;\;} i+j+k=2n-2,\\
0, \mbox{\;\; otherwise.}
\end{cases}\]
\end{Proposition}

\begin{proof} 
Using Equation~\eqref{eq:1-insertion}, we easily deduce that
\[ \widetilde{\langle [\epsilon|\epsilon^i]}\widetilde{\langle [\epsilon|\epsilon^j]}\widetilde{\langle [\epsilon|\epsilon^k]}= \frac{1}{2}\big( \sum_{s+t=k} \delta_{j+s=n-1}\delta_{i+t=n-1}+\delta_{j+t=n-1}\delta_{i+s=n-1}\big)\lambda+O(\lambda^2).\]
Reading off the coefficient of $\lambda$ yields the result.
\end{proof}
\section{Categorical primitive forms}
\label{sec:primitive}

One can obtain more information about the categorical Gromov-Witten invariants by considering a (formal) family of deformations of the potential function $W(x)=\frac{1}{n+1} x^{n+1}$.  Classically one obtains the $g=0$ invariants of singularities by considering a versal family of deformation of the singularity, and constructing a (formal) Frobenius manifold which encodes the $g=0$ invariants.  A key ingredient in this construction is Saito's theory of primitive forms~\cite{Sai2}.

In this section we present a categorical analogue of Saito's definition for $\Z/2\Z$-graded $A_\infty$ algebras.  The results in this section are not specific to categories of matrix factorizations, and apply for any $\Z/2\Z$-graded $A_\infty$ algebra.

\paragraph{\bf Connection in the $u$-direction.} 
Let $A$ be a $\Z/2\Z$-graded strictly unital $A_\infty$ algebra. There is a natural meromorphic $u$-direction connection on $HC^-_*(A)=H_*\big( CC_*(A)\series{u}, b+uB\big)$, the negative cyclic homology of $A$. In general this connection has an order two pole at $u=0$. When the $\Z/2\Z$-grading lifts to a $\Z$-grading it reduces to a simple pole.

To define this connection consider the following deformation of $A$ defined over $t\in \mathbb{A}^1-{0}$:
\[ \mu_k^t(a_1,\cdots,a_k)=t^{2-k}\mu_k(a_1,\cdots,a_k), \forall \; k\geq 0.\]
Its Kodaira-Spencer class is given by
\[ \KS(\partial/\partial t) = \prod_k (2-k) t^{1-k} \mu_k.\]
Its restriction to the fiber $t=1$ defines a class $[\KS(\partial/\partial t)|_{t=1}]\in HH^*(A)$ represented by 
\[ \KS(\partial/\partial t)|_{t=1}= \prod_k (2-k)\mu_k.\] 
Define the $u$-direction connection operator by
\begin{equation}~\label{eq:u-connection}
 \nabla_{\partial/\partial u} := \frac{\partial}{\partial u} +\frac{\Gamma}{2u}+\frac{B^{1|1}(\KS(\partial/\partial t)|_{t=1}; -)}{2u}+\frac{b^{1|1}(\KS(\partial/\partial t)|_{t=1}; -)}{2u^2}.
 \end{equation}
Here the operator $\Gamma$ is defined by $\Gamma(a_0|a_1|\cdots|a_k)=-k a_0|a_1\cdots|a_k$ and the operators $b^{1|1}$ and $B^{1|1}$ are as in~\cite{She}.

\begin{Lemma}~\label{lem:u-connection}
The operator $\nabla_{\partial/\partial u}$ descends to a connection operator
\[ \nabla_{\partial/\partial u}: HC^-_*(A) \ra u^{-2} \cdot HC^-_*(A),\]
with an order two pole at $u=0$. If the $\Z/2\Z$-grading of $A$ extends to a $\Z$-grading then this connection has a simple pole at $u=0$. 
\end{Lemma}

\begin{proof}
Rewrite the operator $2u\nabla_{\partial/\partial u}$ as
\[ 2u\nabla_{\frac{\partial}{\partial u}}=\big( 2u\frac{\partial}{\partial u} +\Gamma+t\frac{\partial}{\partial t}|_{t=1}\big)-\Big( t\frac{\partial}{\partial t}-u^{-1} b^{1|1} (\KS(t\frac{\partial}{\partial t}); -)-B^{1|1}(\KS(t\frac{\partial}{\partial t}); -)\Big)|_{t=1}.\]
Observe that the second term 
\[ \Big( t\frac{\partial}{\partial t}-u^{-1} b^{1|1} (\KS(\partial/\partial t); -)-B^{1|1}(\KS(\partial/\partial t); -)\Big)= \nabla^\GM_{t\frac{\partial}{\partial t}}\]
is the Getzler-Gauss-Manin connection operator~\cite{Get}. Thus, this operator descends to cyclic homology. For the first term one verifies by direct computation that the commutator
\[ [\big( 2u\frac{\partial}{\partial u} +\Gamma+t\frac{\partial}{\partial t}|_{t=1}\big), b+uB]= b+uB.\]
Hence it also descends to cyclic homology. 

Furthermore, if $A$ is $\Z$-graded,  then the $t$-direction deformation is trivialized by the linear $A_\infty$ quasi-isomorphisms $a\mapsto t^{|a|} a$. This implies that the connection operator $\nabla^\GM_{t\frac{\partial}{\partial t}}=0$ in homology, which removes the second order pole at $u=0$.
\end{proof}

\begin{Definition} A choice of a Calabi-Yau structure on an $A_\infty$ algebra $A$ consists of the choice of an element $\omega\in HC_*^-(A)$ such that the induced pairing \[ \langle-,-\rangle: A\otimes A \ra \mathbb{K}, \;\; \langle a, b\rangle := \langle [\mu_2(a,b)], \pi(\omega)\rangle_{{\sf Muk}}\] is non-degenerate. Furthermore, such a class $\omega$ is called primitive if the Hochschild homology class $\pi(\omega)\in HH_*(A)$ generates $HH_*(A)$ as a $HH^*(A)$-module.  We will call such a pair $(A, \omega)$ a {\em primitive Calabi-Yau algebra}.\end{Definition}

\medskip \begin{remark} Another definition of Calabi-Yau structure (of degree $d$) is in terms of an invariant trace map $\Tr: A \rightarrow \mathbb{K}[-d]$. The definition above induces such a trace map when $A$ is minimal. Both definitions appeared in~\cite{KonSoi} already, with the difference that in [loc.\ cit.] the authors used cyclic cohomology instead of cyclic homology. The two approaches are equivalent using the Mukai pairing.  
\end{remark}

\begin{Definition}~\label{def:goodsplitting}
Let $(A,\omega)$ be a primitive Calabi-Yau $A_\infty$ algebra. A splitting $s: HH_*(A) \ra HC_*^-(A)$ of the Hodge filtration is called a {\em good splitting} if the subspace $u^{-1} \imag (s) [u^{-1}]\subset HP_*(A)$ is preserved by the connection $\nabla_{u\frac{\partial}{\partial u}}$. A good splitting $s$ is said to be $\omega$-compatible if there exists a constant $r\in \mathbb{K}$ such that
\[ \nabla_{u\frac{\partial}{\partial u}}\omega - r\omega\in u^{-1} \imag (s).\]
\end{Definition}

\paragraph{{\bf Versal families of $A_\infty$ algebras.}}
Let $R$ be a regular local ring with maximal ideal $\mathfrak{m}$. Let $\cA$ be an $A_\infty$ algebra over $R$. Assume that its underlying $\Z/2\Z$-graded $R$-module is flat and of finite rank over $R$. Furthermore assume that the central fiber $A_\infty$ algebra $A=\cA/\mathfrak{m}\cA$ is minimal. We consider $\cA$ as a formal deformation of $A$.

The Hochschild chain complex of $\cA$ is defined as
\[ CC_*(\cA):= \varprojlim_{N} \big( CC_*(\cA/ \mathfrak{m}^N \cA), b\big),\]
where the differential $b$ is the usual Hochschild differential of the $A_\infty$ algebra $\cA/\mathfrak{m}\cA$. The periodic cyclic chain complex of $\cA$ is defined by the inverse limit \[\varprojlim_{N} \big( CC_*(\cA/ \mathfrak{m}^N \cA)\laurent{u}, b+uB\big) .\] 
Similarly we also define the cyclic chain complex/negative cyclic chain complex as the $\mathfrak{m}$-adic completions of the corresponding chain complexes of $\cA/\mathfrak{m}^N\cA$ .

To summarize, the negative cyclic homology of $\cA$ carries the following structures:
\begin{itemize}
\item[---] The higher residue pairing
\[ \langle-,-\rangle_{\hres}: HC^-_*(\cA) \otimes HC^-_*(\cA) \ra R\series{u}.\]
\item[---] The Getzler-Gauss-Manin connection
\[ \nabla^{\GM}: \Der (R) \otimes HC^-_*(\cA) \ra u^{-1} HC^-_*(\cA).\]
This connection has a first order pole at $u=0$.
\item[---] An extended $u$-direction connection
\[ \nabla_{\frac{\partial}{\partial u}}: HC^-_*(\cA) \ra u^{-2} HC^-_*(\cA).\]
It has a second order pole at $u=0$.
\end{itemize}

\begin{Lemma}
The $u$-direction connection commutes with the Getzler-Gauss-Manin connection in the $R$-direction. That is, for any vector field $v \in \Der (R)$, we have
\[ [\nabla_{\frac{\partial}{\partial u}}, \nabla_v^\GM]=0.\]
\end{Lemma}

\begin{proof}
Since $\nabla_v^\GM$ is $u$-linear it suffices to prove the statement for the operator $2u\nabla_{\partial/\partial u}$. As in the proof of Lemma~\ref{lem:u-connection} we rewrite the operator $2u\nabla_{\partial/\partial u}$ as
\[ 2u\nabla_{\frac{\partial}{\partial u}}=\big( 2u\frac{\partial}{\partial u} +\Gamma+t\frac{\partial}{\partial t}|_{t=1}\big)-\Big( t\frac{\partial}{\partial t}-u^{-1} b^{1|1} (\KS(t\frac{\partial}{\partial t}); -)-B^{1|1}(\KS(t\frac{\partial}{\partial t}); -)\Big)|_{t=1}.\]
The second term 
\[ \Big( t\frac{\partial}{\partial t}-u^{-1} b^{1|1} (\KS(\partial/\partial t); -)-B^{1|1}(\KS(\partial/\partial t); -)\Big)= \nabla^\GM_{t\frac{\partial}{\partial t}}\]
is the Getzler-Gauss-Manin connection operator which commutes with $\nabla^\GM_v$ because the Getzler-Gauss-Manin connection is flat and $[t\frac{\partial}{\partial t}, v]=0$. For the first term commutativity follows from direct computation.
\end{proof}

\begin{Definition}
A formal deformation $\cA$ of $A$ over $R$ is called versal if the Kodaira-Spencer map
\begin{equation*}
 \KS : \Der (R) \ra HH^*(\cA)
 \end{equation*}
is an isomorphism.
\end{Definition}
\bigskip

\textit{Remark.} Using the obstruction theory of $A_\infty$ homomorphisms one can show that a versal deformation of $A$ is unique up to homotopy of $A_\infty$ algebras.  

\begin{Definition}~\label{def:primitive}
Let $(A,\omega)$ be a primitive Calabi-Yau $A_\infty$ algebra. Let $\cA$ be a versal deformation of $A$. An element $\zeta\in HC^-_*(\cA)$ is called a categorical primitive form if it satisfies the following conditions:
\begin{itemize}
\item[P0.] The restriction of $\zeta$ to the central fiber is equal to $\omega$.
\item[P1.] (Primitivity) The map 
\[ \lambda: \Der (R) \ra  HC^-_*(\cA)/u\cdot HC^-_*(\cA)=HH_*(\cA)\]
defined by $v\mapsto \pi\big(u\nabla_v^{\GM}\zeta\big)$  is an isomorphism.
\item[P2.] (Orthogonality) For any tangent vectors $v_1, v_2\in \Der (R)$ we have 
\[\langle u\nabla^{\GM}_{v_1} \zeta, u\nabla^{\GM}_{v_2} \zeta\rangle_{{\hres}} \in R.\]
\item[P3.] (Holonomicity) For any tangent vectors $v_1,v_2,v_3\in \Der (R)$ we have
\[ \langle u\nabla^{\GM}_{v_1} u\nabla^{\GM}_{v_2} \zeta, u\nabla^{\GM}_{v_3} \zeta\rangle_{{\hres}} \in R\oplus u\cdot R.\]
\item[P4.] (Homogeneity) Let $E:=\KS^{-1}(\mu)\in \Der(R)$ be the vector field corresponding to the Hochschild cohomology class represented by the structure map $\mu$ itself, under the Kodaira-Spencer map. Then there exists a constant $r\in \mathbb{K}$ such that
\[ (\nabla_{u\frac{\partial}{\partial u}}+\nabla_E)\zeta= r\zeta.\]
\end{itemize}

\end{Definition}
\medskip

\paragraph{\bf Frobenius Manifolds associated to $(\cA,\zeta)$.} 
Let $\cA$ be a versal deformation of $A$ over a regular local ring $R$, and assume that a categorical primitive form $\zeta$ exists for this family of $A_\infty$ algebras. It is a standard procedure to construct a formal Frobenius manifold structure on the deformation space $\Spec R$. Explicitly, the algebra structure on the tangent space is given by
\[ \frac{\partial}{\partial t_i}\circ \frac{\partial}{\partial t_j}:=\KS^{-1}\big( \KS(\frac{\partial}{\partial t_i})\cup \KS(\frac{\partial}{\partial t_j})\big).\]
Here $\cup$ is the cup product on the Hochschild cohomology. The metric is defined by
\[ \langle \frac{\partial}{\partial t_i}, \frac{\partial}{\partial t_j}\rangle:= \langle u\nabla^\GM_{\frac{\partial}{\partial t_i}}\zeta, u\nabla^\GM_{\frac{\partial}{\partial t_j}}\zeta\rangle_\hres.\]
Furthermore, the vector field $E=\KS^{-1}(\mu)$ as in Definition~\ref{def:primitive} is an Euler vector field. And $e:=\KS^{-1}(\bone)$ is the unit vector field. We refer to Saito-Takahashi~\cite{SaiTak} for details of this construction. One can then define the genus zero potential function $\cF$ by formula
\begin{equation}~\label{eq:potentialfunction}
 \frac{\partial}{\partial \tau_i}\frac{\partial}{\partial \tau_j}\frac{\partial}{\partial \tau_k}\cF=\langle u\nabla^\GM_{\frac{\partial}{\partial \tau_i}} u\nabla^\GM_{\frac{\partial}{\partial \tau_j}}\zeta, u\nabla^\GM_{\frac{\partial}{\partial \tau_k}}\zeta \rangle_\hres \mid_{u=0}.\end{equation}
Here $\tau_i$'s are flat coordinates for the above metric. In general, it is a difficult task to find exact formula of the flat coordinates. In the next section we shall use a perturbative method developed by Li-Li-Saito-Shen~\cite{LLSS} in order to compute these expressions order by order.

\section{Invariants of $(\cA_n,\zeta)$}
\label{sec:computation}

In this section, we prove that for the $A_\infty$ algebra $A_n$ there
exists an essentially unique categorical primitive form. This
coincides with the corresponding statement in the commutative world.
Moreover, the categorical approach is quite computable in this case,
and we prove that the genus zero categorical Gromov-Witten invariants
of $A_n$ match those obtained from Saito's commutative construction of
a Frobenius manifold associated to the function
$f(x)=\frac{1}{n+1}x^{n+1}$.

\paragraph{{\bf Hochschild cohomology.}}
As is well-known, the formal deformation theory of the $A_\infty$
algebra $A_n$ is governed by the shifted Hochschild cochain complex of
$A_n$. Its tangent space is given by the Hochschild cohomology of
$A_n$. Dual to the calculation in Proposition~\ref{prop:homology}, we
have the following result.

\begin{Proposition} The Hochschild cohomology of $A_n$ is an
$n$-dimensional $\Z/2\Z$-graded vector space concentrated in the even
part. It is generated by
\[ HH^*(A_n)= {{\sf span}} ( [\bone], [ \epsilon \mapsto \bone],
  \cdots, [\epsilon^{n-1}\mapsto \bone])\] 
where $[\epsilon^k\mapsto \bone]$ denotes a class represented by a
multi-linear map in $Hom( \overline{A_n}^{\otimes k}, A)$ which sends
$\epsilon^k$ to $\bone$ and all other inputs to zero. 
\end{Proposition}

\paragraph{{\bf Versal deformation of $A_n$.}}
These infinitesimal deformations are easily seen to extend to
formal deformations of $A_n$, as a weakly curved $A_\infty$
algebra. Indeed, denote by
\[ R=\mathbb{K}\series{t_0, \cdots, t_{n-1}}.\]
Define an $A_\infty$ structure on $\cA_n:=A_n\otimes_{\mathbb{K}} R$ by formula
\begin{align}~\label{eq:deformation}
\begin{split}
\mu_0 & = t_0 \bone,\\
\mu_1(\epsilon) &= t_1 \bone,\\
\mu_2(\epsilon,\epsilon) &= t_2 \bone,\\
&\vdots,\\
\mu_{n-1}(\epsilon^{\otimes n-1}) &= t_{n-1} \bone,\\
\mu_{n+1}(\epsilon^{\otimes n+1}) &= \frac{1}{n+1}\bone.
\end{split}
\end{align}
One easily sees that the Kodaira-Spencer class of this family is given by
\[ \KS(\partial/\partial t_j)= [\epsilon^j \mapsto \bone].\] 
We will now study the Hochschild chain complex  of the
family $\cA_n$ in order to compute a categorical primitive form for
this family.

\paragraph{{\bf Rational weight on the Hochschild chain complex.}}
The equations of the versal deformation~\eqref{eq:deformation} imply that if we assign the following rational weights to the variables 
\begin{align*}
\wt(t_j)&:=-j\cdot \frac{2}{n+1}+2, \;\; (j=0,\cdots,n-1),\\
\wt(\epsilon)&:=\frac{n-1}{n+1},
\end{align*}
the $A_\infty$ algebra $\cA_n$ is $\Q$-graded, i.e. the map $\mu_k$
has weight $2-k$. These weights induce rational weights on the
Hochschild chain complex of $\cA_n$ given by
\[ \wt( \bone|\epsilon^k) = \frac{2k}{n+1},\;\; \wt(\epsilon|
  \epsilon^k)= \frac{2k+1-n}{n+1}, \;\;
  \wt(t_j)=j\cdot\frac{2}{n+1}-2.\] 
With respect to these weights, the Hochschild differential has weight
$-1$.

\begin{Proposition}~\label{prop:homologyfamily}
The Hochschild homology of $\cA_n$ is
\[ HH_*(\cA_n)\cong HH_*(A_n)\otimes_{\mathbb{K}} \mathbb{K}\series{t_0,\cdots,t_{n-1}}\]
as a rationally weighted $R$-module. In particular it is a free $R$-module of rank $n$.
\end{Proposition}

\begin{proof}
By Proposition~\ref{prop:homology}, we can pick a homotopy retraction
$(i,p,h): HH_*(A_n) \cong CC_*(A_n)$, which exists since we are over a
field $\mathbb{K}$. We can furthermore choose $i$ and $p$ to have
weight zero, while $h$ has weight $1$. Extend this data
$R/\mathfrak{m}^N$-linearly to get a homotopy retraction which we
still denote by
\[ (i,p,h): HH_*(A_n) \otimes_\mathbb{K} R/\mathfrak{m}^N \cong
CC_*(A_n)\otimes_\mathbb{K} R/\mathfrak{m}^N.\] 
We think of the complex $CC_*(\cA_n/\mathfrak{m}^N\cA_n)$ as a
perturbation of the right hand side. The result follows from
homological perturbation lemma and the fact that $HH_*(A_n)$ is
concentrated at the odd degree, and hence can not support nonzero
differential. 
\end{proof}

\begin{Theorem}
\label{thm:primitive} 
Let $\omega$ be the unique weight-preserving Calabi-Yau $A_\infty$
structure on $A_n$. As in  Proposition~\ref{prop:splitting} $\omega$
is given by
\[ \omega:= s([\epsilon|\epsilon^{n-1}])= \sum_{l\geq 0}^\infty
(-1)^l c_{l} \;\epsilon|\epsilon^{n-1+nl+l} u^l\] 
where $c_{0}=1$ and $c_{l}:=\prod_{j=0}^{l-1} \big((n+j(n+1)\big)$ for any
$l\geq 1$. Let $\cA_n$ be the versal deformation of the algebra $A_n$.
Then there exists a unique categorical primitive form $\zeta\in
HC^-_*(\cA_n)$ lifting $\omega$. 
\end{Theorem}

\begin{proof}
In addition to the above weights assign $\wt(u)=-2$ to the circle
parameter. First, observe that the Euler vector field
$E=\KS^{-1}(\mu)$ is given by 
\[ E= \sum_{j=0}^{n-1} (1-\frac{j}{n+1}) t_j\frac{\partial}{\partial t_j}.\]
Then one verifies that
\[ \nabla_{u\frac{\partial}{\partial u}}+\nabla_E = -\frac{\wt}{2}.\]
This implies that the homogeneity condition P4. is equivalent to the
requirement that the primitive form be homogeneous of weight $-2r$.

Next we prove that there exists a unique weight-preserving primitive form that lift $\omega$. The Calabi-Yau structure $\omega$ of $A_n$ is given by
\[ \epsilon|\epsilon^{n-1}+O(u),\]
which is homogeneous of weight $\frac{n-1}{n+1}$. By the previous proposition, since the weights of the $t_j$ variables are all negative, there is a unique weight-preserving lift of the Hochschild class $[\epsilon|\epsilon^{n-1}]$ to a Hochschild class of $\cA$. Let us denote this Hochschild class by $\Omega_0\in HH_*(\cA_n)$.  Let us consider lifting $\Omega_0$ to a class of the form
\[ \Omega_0+u\Omega_1+u^2\Omega_2+\cdots\]
in the negative cyclic homology of $\cA_n$. The obstruction to the
existence of $\Omega_1$ is a class in $HH_*(\cA)$ represented by
$B\Omega_0$. Since $B$ has weight $1$, and $\Omega_0$ has weight
$\frac{n-1}{n+1}$ this obstruction class has weight
$1+\frac{n-1}{n+1}$. Proposition~\ref{prop:homologyfamily} shows that
this class must be trivial.

To see the uniqueness of the extension $\Omega_1$, observe that the
set of liftings, up to homology, is a torsor over the subgroup of
$HH_*(\cA)$ of classes of weight $2+\frac{n-1}{n+1}$. This group which
vanishes again by Proposition~\ref{prop:homologyfamily}.

The existence and uniqueness of $\Omega_k$ for $k\geq 2$ is similar.
\end{proof}

\textit{Remark.}
The same argument using weights also proves  the uniqueness for all
ADE type singularities. 

\paragraph{{\bf Comparison with Saito's genus zero invariants.}} 
In the remainder of this section we prove that the genus zero
categorical Gromov-Witten invariants of $\MF(\frac{1}{n+1}x^{n+1})$
obtained from the above weight-preserving categorical primitive form
match with Saito's invariants of $\frac{1}{n+1}x^{n+1}$. To obtain
this result we apply a reconstruction theorem of
Fan-Jarvis-Ruan~\cite{FJRW} which asserts that the genus zero Saito
invariants can be recovered from the two-point and three-point
functions along with the single four-point function
$\langle x,x,x^{n-1},x^{n-1}\rangle$. 

To obtain this reconstruction result Fan-Jarvis-Ruan use two formal
properties of the genus zero Gromov-Witten potential: the WDVV
equation and the dimension axiom, both of which are known to hold for
Saito's invariants.

In general it is well-known that WDVV equation is equivalent  to the
associativity condition in the definition of Frobenius manifold.  Thus
the WDVV equation also holds in the categorical case.  The dimension
axiom for the categorical invariants can be deduced from the existence
of an Euler vector field. This is proved in detail in~\cite{SaiTak}.
In the $A_n$ case, this gives the equation 
\[ E(\cF)=(3- \frac{n-1}{n+1})\cF\] 
where $\cF$ stands for the genus zero potential defined using
Equation~\ref{eq:potentialfunction}.

Thus, both the categorical and Saito's invariants are determined by the  two-point, three-point, and a particular four-point function.  For Saito's invariants with primitive
form $dx$ these quantities are as follows, see~\cite{JKV}:
\begin{itemize} \item[(1)]
  $\langle x^i, x^j \rangle=\delta_{i+j=n-1}$, \item[(2)]
  $\langle x^i, x^j, x^k\rangle=\delta_{i+j+k=n-1}$, \item[(3)]
  $\langle x, x, x^{n-1}, x^{n-1}\rangle=1$. \end{itemize} 
Using the results of Section~\ref{sec:invariants} these calculations
match the categorical two-point and
three-point functions (Propositions~\ref{prop:pairing} and~\ref{prop:coproduct}) under the isomorphism
\[ \Phi: HH_{*+1}(A_n) \cong \mathbb{K}[x]/(x^n), \;\;
  \Phi(\epsilon|\epsilon^i)=x^{n-1-i}.\]

We will use a method developed by Li-Li-Saito-Shen~\cite{LLSS} to
compute the categorical four-point function 
\[ \langle \Phi^{-1}(x), \Phi^{-1}(x), \Phi^{-1}(x^{n-1}),
  \Phi^{-1}(x^{n-1})\rangle. \]
Their method relies on constructing a trivialization of the
Gauss-Manin connection on the twisted de Rham cohomology. The
following proposition is a non-commutative analogue, giving a
trivialization of the Getzler-Gauss-Manin connection on $HP_*(\cA_n)$.
To set it up, for $0\leq i \leq n-1$, denote by $\phi_i\in C^*(\cA_n)$
the Hochschild cochain $[\epsilon^k\mapsto \bone]$ of $A_n$ extended
$R$-linearly to $\cA_n$. To simplify the notation denote by 
\[ b_i:=b^{1|1}(\phi_i), \;\; 0\leq i\leq n-1\] 
the action of Hochschild cochain on Hochschild chains. The proof of
the following proposition is a direct computation which we shall omit.

\begin{Proposition}~\label{prop:trivialization}
We have an isomorphism of differential modules over $R$:
\[ \exp\big(-\frac{\sum_{i=0}^{n-1} t_i b_i}{u}\big): \big(HP_*(\cA_n), \nabla^\GM\big) \ra \big( HP_*(A_n)\otimes_\mathbb{K} R, d=\sum_{j=0}^{n-1}\frac{\partial}{\partial t_j} dt_j\big).\]
\end{Proposition}

\begin{Theorem}
\label{thm:comparison} 
Consider the Calabi-Yau $A_\infty$ category
$\MF(\frac{1}{n+1}x^{n+1})$ and its versal deformation $\cA_n$. Endow
this family with the (unique up to scalars) categorical primitive form
$\zeta$ described in Theorem~\ref{thm:primitive}. Then the categorical
genus zero invariants computed from $\cA_n, \zeta)$ agree with the
corresponding invariants obtained through Saito's procedure of
unfolding of singularities of $\frac{1}{n+1}x^{n+1}$, using the
(unique up to scalar) primitive form $dx$. 
\end{Theorem}

\begin{proof}
From the above discussion, the result will follow if we succeed in
matching the categorical Gromov-Witten four-point invariant with
Saito's invariant
\[ \langle x, x, x^{n-1}, x^{n-1}\rangle=1. \]
Li-Li-Saito-Shen deduced an inductive formula to read off
Gromov-Witten invariants from a primitive form $\zeta$ by solving the
equation 
\begin{equation}~\label{eq:formula}
 \exp\big(-\frac{\sum_{i=0}^{n-1} t_i b_i}{u}\big) \zeta = \omega + u^{-1} J_{-1} +u^{-2} J_{-2} +\cdots, \;\; J_{-k} \in \bigoplus_{j=0}^{n-1} R\cdot s_j
\end{equation}
order by order, with respect to the $\mathfrak{m}$-adic filtration.
Here the $s_j$'s are the splittings defined in
Proposition~\ref{prop:splitting}. We write
$\zeta=\zeta^{(0)}+\zeta^{(1)}+\cdots$ so that $\zeta^{(k)}$ is the
order $k$ part of $\zeta$. Since the primitive form restricts to
$\omega=s_{n-1}$ on the central fiber we have $\zeta^{(0)}=s_{n-1}$.

To first order the above equation yields
\[ \zeta^{(1)} -\frac{\sum_{i=0}^{n-1} t_i b_i}{u} s_{n-1}= u^{-1} J_{-1}^{(1)}.\]
This equation is uniquely solved by setting
\[ J_{-1}^{(1)}= - \sum_{i=0}^{n-1} t_i s_{n-1-i}, \mbox{\; and \;} \zeta^{(1)} = u^{-1} \sum_{i=0}^{n-1} t_i \big( b_i (s_{n-1}) -s_{n-1-i}\big).\]

In order to get the four-point function we need to compute $J_{-1}$ up
to order $2$, and $J_{-2}$ up to order $3$. These can be derived again
using Equation~\eqref{eq:formula}. To second order it gives
\[ \zeta^{(2)} - \frac{\sum_{i=0}^{n-1} t_i b_i}{u} \zeta^{(1)} + \frac{1}{2}\frac{\sum_{i=0}^{n-1} t_i b_i}{u}\frac{\sum_{j=0}^{n-1} t_j b_j}{u} s_{n-1} = u^{-1} J_{-1}^{(2)} + u^{-2} J_{-2}^{(2)}.\]
From this, and using formula of $s_j$'s in Proposition~\ref{prop:splitting}, we deduce that
\begin{align*}
J_{-2}^{(2)} & = \frac{1}{2} \sum_{0\leq i,j\leq n-1 \atop i+j\leq n-1} t_it_j \cdot s_{n-1-i-j}\\
J_{-1}^{(2)} &= \sum_{0\leq i,j\leq n-1 \atop i+j\geq n+1} \big( j-\frac{n}{2}\big) t_it_j \cdot s_{2n-i-j}.
\end{align*}

To third order Equation~\eqref{eq:formula} is
\begin{align*} 
&\zeta^{(3)} - \frac{\sum_{i=0}^{n-1} t_i b_i}{u} \zeta^{(2)} + \frac{1}{2}\big(\frac{\sum_{i=0}^{n-1} t_i b_i}{u}\big)^2\zeta^{(1)}-\frac{1}{6}\big(\frac{\sum_{i=0}^{n-1} t_i b_i}{u}\big)^3 s_{n-1}  \\
&= u^{-1} J_{-1}^{(3)} + u^{-2} J_{-2}^{(3)}+u^{-3} J_{-3}^{(3)}.
\end{align*}
We only extract $J_{-2}^{(3)}$ from it:
\[ J_{-2}^{(3)}= \sum_{0\leq i,j,k\leq n-1 \atop i+j+k\geq n+1} \big( \frac{n}{6}-\frac{k}{2}\big)t_it_jt_k\cdot s_{2n-i-j-k}.\]

The key point of this order by order algorithm is the fact that the
term $J_{-1}$ gives flat coordinates of the Frobenius manifold
associated to $\zeta$.
\begin{align}~\label{eq:coordinates}
\begin{split}
\tau_k &= -t_k+ \sum_{k+2\leq j\leq n-1} (j-\frac{n}{2})\cdot t_{n+1+k-j} t_j+O(t^3), \;\; 0\leq k\leq n-3,\\
\tau_{n-2} &=-t_{n-2}+O(t^3),\\
\tau_{n-1} &=-t_{n-1}+O(t^3).
\end{split}
\end{align}
That is, in flat coordinates, we have
$J_{-1}= \sum_{0\leq i\leq n-1} \tau_i\cdot s_{n-1-i}$.

Let us write the term $J_{-2}$ as \[ J_{-2} = \sum_l J_{-2}^l s_l.\]
Then the coefficients $J_{-2}^l$ give the partial derivatives of the
genus zero potential $\cF$ in flat coordinates. More precisely, we
have
\[ \frac{\partial}{\partial \tau_l} \cF = \sum_{0\leq k\leq n-1}
  g_{kl} J_{-2}^k(\tau).\] Here the constants $g_{kl}$ are given by
$g_{kl}=\langle [\epsilon|\epsilon^{k}],
[\epsilon|\epsilon^{n-1-l}]\rangle_{{\sf Muk}}= \delta_{k=l}$ by
Proposition~\ref{prop:pairing}. (The discrepancy in $g_{kl}$ is due to
the reverse indices in the formula of $J_{-1}$.) This implies that
\[ \frac{\partial}{\partial\tau_l} \cF = J_{-2}^{l}(\tau).\] Using the
formulas for $J_{-2}^{(2)}$, $J_{-2}^{(3)}$ we obtain that
\begin{align*} \frac{\partial}{\partial\tau_l} \cF &= \frac{1}{2}
  \sum_{0\leq i\leq n-1-l} t_it_{n-1-l-i} +\sum_{0\leq i,j,k\leq
    n-1\atop i+j+k=2n-l} (\frac{n}{6}-\frac{k}{2}) t_it_jt_k+O(t^4).
\end{align*} Using the above coordinate change
formula~\eqref{eq:coordinates} and setting $l=n-1$ we obtain
\[ \frac{\partial}{\partial\tau_{n-1}} \cF = \frac{1}{2}
  \tau_0^2+\sum_{2\leq j\leq n-1}(n-2j) \tau_0
  \tau_j\tau_{n+1-j}-\sum_{i+j+k=n+1}
  (\frac{n}{6}-\frac{k}{2})\tau_i\tau_j\tau_k +O(\tau^4).\] The
four-point function
$\frac{\partial}{\partial\tau_1}\frac{\partial}{\partial\tau_1}\frac{\partial}{\partial\tau_{n-1}}\frac{\partial}{\partial\tau_{n-1}}
\cF |_{\tau=0}$ is then given by
\[
  \frac{\partial}{\partial\tau_1}\frac{\partial}{\partial\tau_1}\frac{\partial}{\partial\tau_{n-1}}\frac{\partial}{\partial\tau_{n-1}}
  \cF |_{\tau=0}=1.\] \end{proof}

\end{document}